\documentclass[final]{elsarticle}
\usepackage{latexsym,bm,amssymb,amsfonts,amsbsy,amsmath,graphics,graphicx,color,url}
\def\dd{\,\mathrm{d}}

\newtheorem{remark}{Remark}
\newtheorem{lemma}{Lemma}

\def \RR {{\mathbb{R}}}

\def \P {\hat{P}}

\def \pmatrix{ \left( \begin{array} }
\def \endpmatrix{ \end{array} \right) }

\def\II{{\cal I}}
\def\PP{{\cal P}}
\def\bfq{{\bf q}}
\def\bfp{{\bf p}}
\def\bfy{{\bf y}}
\def\bgam{\bm{\gamma}}
\def\bfeta{\bm{\eta}}
\def\diag{\mbox{diag}}

\begin{document}

\title{A note on the efficient implementation of Hamiltonian BVMs}
\tnotetext[t1]{Work developed within the project ``Numerical methods and
software for differential equations''.}

\author[lb]{Luigi Brugnano}\ead{luigi.brugnano@unifi.it}
\author[fi]{Felice Iavernaro}\ead{felix@dm.uniba.it}
\author[dt]{Donato Trigiante}\ead{trigiant@unifi.it}

\address[lb]{Dipartimento di Matematica,
Universit\`a di Firenze, Viale Morgagni 67/A, 50134 Firenze
(Italy).}
\address[fi]{Dipartimento di Matematica, Universit\`a di
Bari, Via Orabona  4,  70125 Bari (Italy).}
\address[dt]{Dipartimento di Energetica, Universit\`a di
Firenze, Via Lombroso 6/17, 50134 Firenze (Italy).}

\begin{abstract} We discuss the efficient implementation of Hamiltonian BVMs
(HBVMs), a recently introduced class of energy preserving methods for canonical
Hamiltonian systems (see \cite{BIT0} and references therein), via their {\em
blended} formulation. We also discuss the case of separable problems, for
which the structure of the problem can be exploited to gain efficiency.
\end{abstract}


\maketitle

\section{Introduction} The conservation of energy allows to avoid the numerical
drift observed when using standard numerical methods for solving canonical
Hamiltonian problems, i.e., problems in the form
\begin{equation}\label{ham}
  y'= J \nabla H(y), \qquad
  J=\pmatrix{cc}0&I_m\\-I_m&0\endpmatrix, \qquad
y(t_0)=y_0\in\RR^{2m},\vspace{-.5em}
\end{equation}
where $H(y)$ is a smooth scalar function and, in general, $I_r$
will hereafter denote the identity matrix of dimension $r$ (when
the lower index will be omitted, the size of the matrix can be
deduced from the context). In this respect, {\em Hamiltonian
Boundary Value Methods (HBVMs)} is a recently introduced class of
methods able to conserve energy when $H(y)$ is a polynomial of
arbitrarily high degree. Clearly, this implies a {\em practical}
conservation of energy for any suitably regular Hamiltonian
function, which will be assumed hereafter. We refer to
\cite{BIT0,BIT1,BIT2,BIT3,BIT4,BIT5}, and references therein, for
an overview on energy-conserving methods and the derivation of
HBVMs. When problem (\ref{ham}) is separable, i.e., when
\begin{equation}\label{sepa}
H(y) \equiv H(q,p) = \frac{1}2 p^Tp - U(q),\qquad
q,p\in\RR^m,\end{equation} then (\ref{ham}) reduces to a special
second order equation,
$$q'' = \nabla U(q),$$ and the associated HBVM may be properly formulated in
order to take advantage, in terms of efficiency, from the above
simplification.

In this paper we investigate the efficient implementation of
HBVMs, also in the case of separable problems. In more details, in
Section~\ref{hbvm} we briefly derive HBVMs. Then, in
Section~\ref{blend} we investigate the efficient solution of the
generated discrete problem, via the {\em blended} implementation
of the methods, which has already proved to be very effective in
other settings (see, e.g.,
\cite{B00,BM02,BM04,BM07,BM09,BM09a,BMM06,BT01}). The case of
separable problems is then discussed in Section~\ref{seprob}. A
few numerical tests, along with some concluding remarks are then
given in Section~\ref{numtest}.

\section{Hamiltonian BVMs (HBVMs)}\label{hbvm}
The derivation of HBVMs will be done according to the approach
followed in \cite{BIT4,BIT5}, which further simplifies the already
simple idea initially used in \cite{BIT0,BIT1,BIT2,BIT3} (see also
\cite{IP08,IT09}). Let us then consider the restriction of problem
(\ref{ham}) to the interval $[t_0,t_0+h]$, with the right-hand
side expanded along an orthonormal basis $\{\P_j\}_{j\ge0}$:
\begin{equation}\label{series}y'(t_0+\tau h) = J\sum_{j\ge0} \P_j(\tau) \int_0^1
\P_j(c)\nabla H(y(t_0+ch))\dd c, \qquad
\tau\in[0,1].\end{equation} In particular, we here consider an
orthonormal polynomial basis, provided by the shifted and scaled
Legendre polynomials on the interval $[0,1]$, even though the
arguments can be easily extended to  more general bases. The basic
idea, is now that of looking for an approximate solution belonging
to the set of polynomials of degree not larger than $s$. This is
achieved by truncating the series at the right-hand side in
(\ref{series}), thus obtaining the approximate problem
\begin{equation}\label{trunc}
\sigma'(t_0+\tau h) = J\sum_{j=0}^{s-1} \P_j(\tau) \int_0^1 \P_j(c)\nabla
H(\sigma(t_0+ch))\dd c, \qquad \tau\in[0,1], \qquad
\sigma(t_0)=y_0.\end{equation}
The approximation to $y(t_0+h)$ is then given by
\begin{equation}\label{y1}y_1\equiv \sigma(t_0+h).\end{equation}The
method can be easily seen to be energy-preserving since, considering that $J$ is
skew-symmetric,
\begin{eqnarray*}
H(y_1) - H(y_0)&=&h\int_0^1 \nabla H(\sigma(t_0+\tau
h))^T\sigma'(t_0+\tau h)\dd\tau\\&=&h\sum_{j=0}^{s-1} \left[
\int_0^1 \P_j(\tau)\nabla H(\sigma(t_0+\tau h))\dd
\tau\right]^TJ\left[ \int_0^1 \P_j(c)\nabla H(\sigma(t_0+ch))\dd
c\right]=0.\end{eqnarray*}

Integrating both sides of the first equation in \eqref{trunc}
yields
\begin{equation}
\label{trunc_int}
 \sigma(t_0+\tau h) = y_0 + h\sum_{j=0}^{s-1}
\int_0^{\tau}\P_j(x)\dd x \int_0^1 \P_j(c)J \nabla
H(\sigma(t_0+ch))\dd c,
\end{equation}
which may be exploited to determine the shape of the unknown
polynomial $\sigma$, provided that a technique to handle the
rightmost integrals is taken into account: the obvious choice is
the use of quadrature formulae. If we assume that $H(y)$ is a
polynomial of degree $\nu$, then the integrals appearing in
(\ref{trunc}) can be exactly computed by a Gaussian formula with
$k$ abscissas $\{c_i\}$, in the event that
\begin{equation}\label{klargenough}k\ge s\nu/2,\end{equation} thus
obtaining a discrete problem in the form
\begin{equation}\label{discr}
 \sigma(t_0+c_ih)\equiv \sigma_i = y_0 + h\sum_{j=0}^{s-1}
\int_0^{c_i}\P_j(x)\dd x \sum_{\ell=1}^k b_\ell \P_j(c_\ell)J\nabla
H(\sigma_\ell), \qquad i=1,\dots,k,
\end{equation}
where the $\{b_i\}$ are the quadrature weights of the formula
defined over the abscissae $\{c_i\}$. For general, suitably
regular (e.g., analytical) Hamiltonian functions, we can still use
formula \eqref{discr} in place of \eqref{trunc_int}, provided that
the integrals in \eqref{trunc_int} are approximated to machine
precision\footnote{As we will see, increasing the order of the
quadrature formula, namely the integer $k$, will not result into
an increase of the computational cost associated with the
implementation of the method.}: in the following, we will always
assume such an accuracy level when a non polynomial function is
considered, and consequently we will make no distinction between
the integrals and the corresponding approximations as well as
between the two polynomials $\sigma$ obtained by solving either
\eqref{discr} or \eqref{trunc_int} (see \cite{BIT5} for more
details).

Method \eqref{discr}-\eqref{y1} is called HBVM($k$,$s$): it was
shown \cite{BIT2} that its  order is $2s$, for all $k\ge s$. In
particular, for $k=s$ it reduces to the well known $s$-stages Gauss
method.

By introducing the matrices ~$\Omega = \diag(b_1,\dots,b_k)$~ and
$$\II_{s-1} =\left( \int_0^{c_i} \P_{j-1}(x)\dd
x\right)_{\scriptsize\begin{array}{l}i=1\dots k\\ j=1\dots
s\end{array}}\in\RR^{k\times s}, \qquad \PP_{r-1} = \left(
\P_{j-1}(c_i)\right)_{\scriptsize\begin{array}{l}i=1\dots k\\
j=1\dots r\end{array}}\in\RR^{k\times r},$$ the HBVM($k$,$s$)  can
be recast as a Runge-Kutta method with Butcher tableau
\begin{equation}\label{buttab}
\begin{array}{c|c} \begin{array}{c} c_1\\ \vdots \\ c_k\end{array} &
A\equiv\II_{s-1}\PP_{s-1}^T\Omega \\ \hline & \begin{array}{ccc}b_1
&\cdots&b_k\end{array}\end{array}
\end{equation}

The next result follows from well-known properties of Legendre
polynomials (hereafter $e_i$ denotes the $i$th unit vector in
$\RR^s$).
\begin{lemma}\label{lem1}\begin{equation}\label{Is}
\II_{s-1} = \PP_s \hat{X}_s\equiv\PP_s\pmatrix{c} X_s\\ \xi_s
e_s^T\endpmatrix,\end{equation}
where \begin{equation}\label{Xs}X_s = \pmatrix{rrrr}
\frac{1}2 & -\xi_1\\
\xi_1     &0   &\ddots\\
          &\ddots &\ddots &-\xi_{s-1}\\
          &&\xi_{s-1} &0\endpmatrix,\qquad \xi_i = \frac{1}{2\sqrt{4j^2-1}},
\quad i\ge1.\end{equation}
\end{lemma}
Consequently, the matrix in the Butcher tableau (\ref{buttab}) can
be written as \begin{equation}\label{A}A = \PP_s\hat
X_s\PP_{s-1}^T\Omega.\end{equation}

Notice that, since $\PP_s \hat X_s$ has $s$ linearly independent
columns, the $k\times k$ coefficient matrix $A$ has rank $s$: it is
then possible to recast the discrete problem in a more convenient
form, which clearly shows that its (block) size is $s$, rather than
$k$ (see also \cite{BIT1}). For this purpose, let us define the
(block) vectors (see (\ref{trunc}) and (\ref{discr}))
\begin{equation}
\label{gammaj}
\bfy = \pmatrix{c} \sigma_1 \\
\vdots \\ \sigma_k\endpmatrix, \quad \bgam = \pmatrix{c} \gamma_0\\ \vdots \\
\gamma_{s-1}\endpmatrix, \quad \gamma_j = \sum_{\ell=1}^k
b_\ell\P_j(c_\ell) J \nabla H(\sigma(t_0+c_\ell h)), \quad
j=0,\dots,s-1.
\end{equation}

 In view of (\ref{trunc}), we see that the vectors  $\gamma_j$
 may be interpreted as the coefficients in the expansion of the
 degree $s-1$ polynomial $\sigma'(t_0+\tau h)$ along the orthonormal
 basis $\{\hat P_j\}_{j=0,\dots,s-1}$.

 From \eqref{discr} one obtains
\begin{equation}\label{ygam}\bfy = e\otimes y_0 +h\II_{s-1}\otimes
I_{2m}\, \bgam,\end{equation} with $e=(1,\dots,1)^T\in\RR^k$, and
then, by virtue of (\ref{gammaj}), one has to solve the equation in the unknown $\bgam$
\begin{equation}\label{eqgam}F(\bgam) \equiv \bgam-
\left(\PP_{s-1}^T\Omega \otimes J\right) \nabla H\left( e\otimes y_0
+ h\II_{s-1} \otimes I_{2m}\bgam\right) = 0.\end{equation} The
application of the simplified Newton iteration for solving
(\ref{eqgam}) yields\begin{equation}\label{it0} \left[ I
-h\PP_{s-1}^T\Omega \II_{s-1} \otimes G_0\right] \Delta^\ell =
-F(\bgam^\ell), \qquad \bgam^{\ell+1} =
\bgam^\ell+\Delta^\ell,\end{equation}  with ~$G_0=\left(J \nabla^2
H(y_0) \right)$.~ By virtue of (\ref{Is}), and considering that
\begin{equation}\label{orto}\PP_{s-1}^T\Omega
\PP_s = \left(I_s~0\right)\in\RR^{s\times s+1},\end{equation} (\ref{it0})
reduces to
\begin{equation}\label{itgam}
\left[ I - hX_s\otimes G_0\right]
\Delta^\ell = -F(\bgam^\ell), \qquad
\bgam^{\ell+1}=\bgam^\ell+\Delta^\ell, \qquad
\ell=0,1,\dots,\end{equation} which, as is readily seen, has
(block) size $s$, rather than $k$.

\section{Blended implementation}\label{blend}
 From the arguments in the previous section, one then concludes that
the discrete problem, to be solved at each integration step when
approximating the Hamiltonian problem (\ref{ham}), is given by
(\ref{eqgam}), thus requiring the solution of (\ref{itgam}). We
are going to solve such equation by means of a {\em blended}
implementation of the method, according to
\cite{B00,BM02,BM04,BT01}. Indeed, such implementation of block
implicit methods has proved to be very effective, leading to the
development of the codes {\tt BiM} \cite{BM04} and {\tt BiMD}
\cite{BMM06} for stiff ODE IVPs and linearly implicit DAEs (the
codes are available at the url \cite{BIMURL}). Let us, for sake of
simplicity, discard the iteration index. Consequently, we have to
solve the linear system
\begin{equation}\label{eq1}\left( I-h X_s \otimes G_0 \right) \Delta =
-F(\bgam) \equiv \bfeta.
\end{equation} Considering that matrix $X_s$ (see (\ref{Xs})) is nonsingular,
such equation can be equivalently written as
\begin{equation}\label{eq2}\rho\left( X_s^{-1}\otimes
I_{2m} -hI_s\otimes G_0\right) \Delta = \rho X_s^{-1}\otimes
I_{2m}\, \bfeta \equiv \bfeta_1,\end{equation} where $\rho$ is a
positive constant. By introducing the (matrix) weight function
\begin{equation}\label{tet}\theta = I_s\otimes \Sigma_0^{-1}, \qquad
\Sigma_0 = (I_{2m}-\rho h G_0)^{-1},
\end{equation} we then obtain the following problem, which has still the same
solution as (\ref{eq1}):
\begin{equation}\label{blendeq} T(\Delta) ~\equiv~
 \theta \left[\left(I-h X_s\otimes G_0\right)\Delta -\bfeta\right] +
(I-\theta)\left[\rho\left( X_s^{-1}\otimes I_{2m} -h I_s\otimes
G_0\right)  \Delta -\bfeta_1\right]~=~0.\end{equation}
One easily realizes that it is obtained as the {\em blending}, with weights
$\theta$ and $(I-\theta)$, of the two equivalent problems (\ref{eq1}) and
(\ref{eq2}), respectively. Problem (\ref{blendeq}) defines the
{\em blended method} associated with the original one, which we
call {\em blended HBVM}, in the present case. The free parameter
$\rho$ is chosen in order to optimize the convergence properties
of the corresponding {\em blended iteration},
\begin{equation}\label{blit}
\Delta_{n+1} = \Delta_n -\theta T(\Delta_n), \qquad
n\ge0,
\end{equation}
with an obvious meaning of the lower index. Such iteration only
requires (see (\ref{tet})) the factorization of the matrix $\Sigma_0$ having
the same size as that of the continuous problem. According to the
linear analysis of convergence in \cite{BM09}, the free parameter
$\rho$ is chosen as
\begin{equation}\label{roopt} \rho = \rho_s \equiv \min\{
|\lambda|\,: ~ \lambda\in\sigma(X_s)\},\end{equation} which
provides optimal convergence properties (in particular, an $L$-convergent
iteration \cite{BM09}). A few values of (\ref{roopt}) are listed in the table
below, for sake of completeness.

\medskip
\centerline{\begin{tabular}{|c|cccc|}
\hline
$s$ & $2$ & $3$ & $4$ & $5$ \\
\hline
$\rho_s$ & $0.2887$ & $0.1967$ & $0.1475$ & $0.1173$ \\
\hline
\end{tabular}}

\begin{remark}
 A nonlinear version of (\ref{blit}) can be readily derived, by
taking $\Delta_n=0$ and updating the vectors $\bfeta$ and
$\bfeta_1$ in (\ref{blendeq}) at each iteration.
\end{remark}

\section{The case of separable problems}\label{seprob}
Let us now apply the method to the
separable problem (\ref{sepa}). By setting the (block) vectors
$$\bfq = \pmatrix{ccc} q_1^T, &\dots~, &q_k^T\endpmatrix^T, \qquad \bfp =
\pmatrix{ccc} p_1^T, &\dots~, &p_k^T\endpmatrix^T,$$one then
obtains (see (\ref{A})),
$$\bfq = e\otimes q_0 +hA\otimes I_m\,\bfp, \qquad \bfp = e\otimes p_0 +hA\otimes
I_m\,\nabla U(\bfq),$$ i.e., since $Ae=c \equiv
(c_1,\dots,c_k)^T$,
 \begin{equation}\label{qprob}\bfq = e\otimes q_0 +h c\otimes p_0 +h^2
A^2\otimes I_m\, \nabla U(\bfq).\end{equation} Moreover, taking
into account (\ref{buttab})--(\ref{A}) and (\ref{orto}), one obtains
\begin{equation}\label{A2}
A^2= \II_{s-1} X_s \PP_{s-1}^T\Omega\,.\end{equation} The new
approximations to $q(t_0+h)$ and $p(t_0+h)$ are then given by
$$q_0 + h p_0 +h^2 e^T\Omega A\otimes I_m\, \nabla U(\bfq), \qquad p_0+h
e^T\Omega\otimes I_m\,\nabla U(\bfq),$$ respectively. By using
similar arguments as those given in Section~\ref{hbvm} (see
(\ref{ygam})), we set $$\bfq = e\otimes q_0 +hc\otimes p_0 +
h^2\II_{s-1}X_s\otimes I_m\, \bgam,$$ thus obtaining the following
equation (which is the analogous of (\ref{eqgam})):
\begin{equation}\label{eqgam2}F(\bgam)
\equiv \bgam -\left(\PP_{s-1}^T\Omega \otimes I_m\right) \nabla
U\left( e\otimes q_0 + hc\otimes p_0 + h^2\II_{s-1}X_s\otimes I_m
\bgam\right) = 0.\end{equation} Similarly as what seen in
Section~\ref{blend}, the application of the simplified Newton
iteration for solving (\ref{eqgam2}) then gives, by virtue of
(\ref{Is}) and (\ref{orto}), and setting $G_0=\nabla^2U(q_0)$,
\begin{equation}\label{itgam2}
\left[ I - h^2X_s^2\otimes G_0\right] \Delta^\ell =
-F(\bgam^\ell), \qquad \bgam^{\ell+1}=\bgam^\ell+\Delta^\ell,
\qquad \ell=0,1,\dots,\end{equation} which, as in the previous
case, has (block) size $s$, rather than $k$. The problem is then exactly that
seen in (\ref{itgam}), via the formal substitutions
\begin{equation}\label{formal}h\longrightarrow h^2, \qquad
X_s\longrightarrow X_s^2.\end{equation} This means that we can
repeat similar steps for the {\em blended} solution of
(\ref{itgam2}), by following the same arguments seen in
Section~\ref{blend}. In more details, (\ref{eq1})--(\ref{blit})
can be repeated, by considering the formal substitutions
(\ref{formal}) and, moreover, $$\rho\longrightarrow\rho^2, \qquad
I_{2m}\longrightarrow I_m.$$ Also in this case \cite{BM07,BM09},
the optimal choice of the parameter $\rho$ turns out to be given
by (\ref{roopt}).

\section{Numerical Tests}\label{numtest}
We here consider a model problem to test the proposed algorithms
and methods, in order to confirm the usefulness of the proposed
approach. In particular, it is clear that a Newton-type iteration,
like (\ref{itgam}) and (\ref{itgam2}), works well when the linear
part of the problem is significant. For this purpose, we consider
the following polynomial Hamiltonian,
\begin{equation}\label{Hqp}
 H(q,p) = \frac{1}2p^2 -10^4 q^2 \left( \frac{4}5 q^3 -\frac{3}4 q^2
-\frac{2}3q +\frac{1}2 \right),
\end{equation} from which we derive the following special second order problem,
\begin{equation}\label{ex}
 q'' = 10^4q\left( 4q^3 -3q^2 -2q + 1\right), \qquad t\in[0,100], \qquad
q(0)=0,\qquad q'(0) = 1.
\end{equation}
For solving (\ref{ex}), we use the following fourth-order numerical methods:
\begin{itemize}
\item the symplectic 2-stages Gauss method (GAUSS2);

\item HBVM(8,2) which is energy conserving, for the problem at hand.
\end{itemize}
For both methods, we consider a fixed-step implementation with stepsize $h$,
with the generated discrete problems solved either with a fixed-point iteration
or with a blended iteration, which  have approximately the same cost, in terms
of function evaluations. Moreover, we also compare the second order
implementation described in Section~\ref{seprob}, with the equivalent first
order Hamiltonian formulation of the problem, as described in
Section~\ref{blend}. Table~\ref{tab1} summarizes the obtained results, in terms
of total number of iterations (blended or fixed-point) for covering the
specified integration interval. From the listed results, one deduces that the
second order formulation of the problem is less demanding in terms of needed
iterations. Moreover, the blended iteration turns out to be both more efficient
and robust than the fixed-point iteration.

\begin{table}[t]
\caption{\protect\label{tab1} total number of iterations for solving the
discrete problems with the specified stepsize $h$  (--  if no convergence).}
{\small
 \begin{tabular}{|c|r|r|r|r|r|r|r|r|}
\hline
 & \multicolumn{4}{c|}{GAUSS2} & \multicolumn{4}{c|}{HBVM(8,2)}\\
 $h$   & \multicolumn{2}{c|}{second order} & \multicolumn{2}{c|}{first
order} & \multicolumn{2}{c|}{second order} & \multicolumn{2}{c|}{first order}\\
& blended & fixed-point & blended & fixed-point & blended & fixed-point &
blended & fixed-point \\
\hline $10^{-3}$ & 664545 &690197 & 952902 & 1217673 &660317&
695765 & 947618 &1225318
 \\
$5\cdot 10^{-3}$ & 242844 & --  & 308406 & -- & 228242 & 223883 &
293949 & 424402\\
$10^{-2}$ & -- & -- & -- & -- & 194163 & -- & 253049 & -- \\
\hline
\end{tabular}
}
\end{table}

\begin{figure}[ht]
\centerline{
\includegraphics[width=7cm,height=7cm]{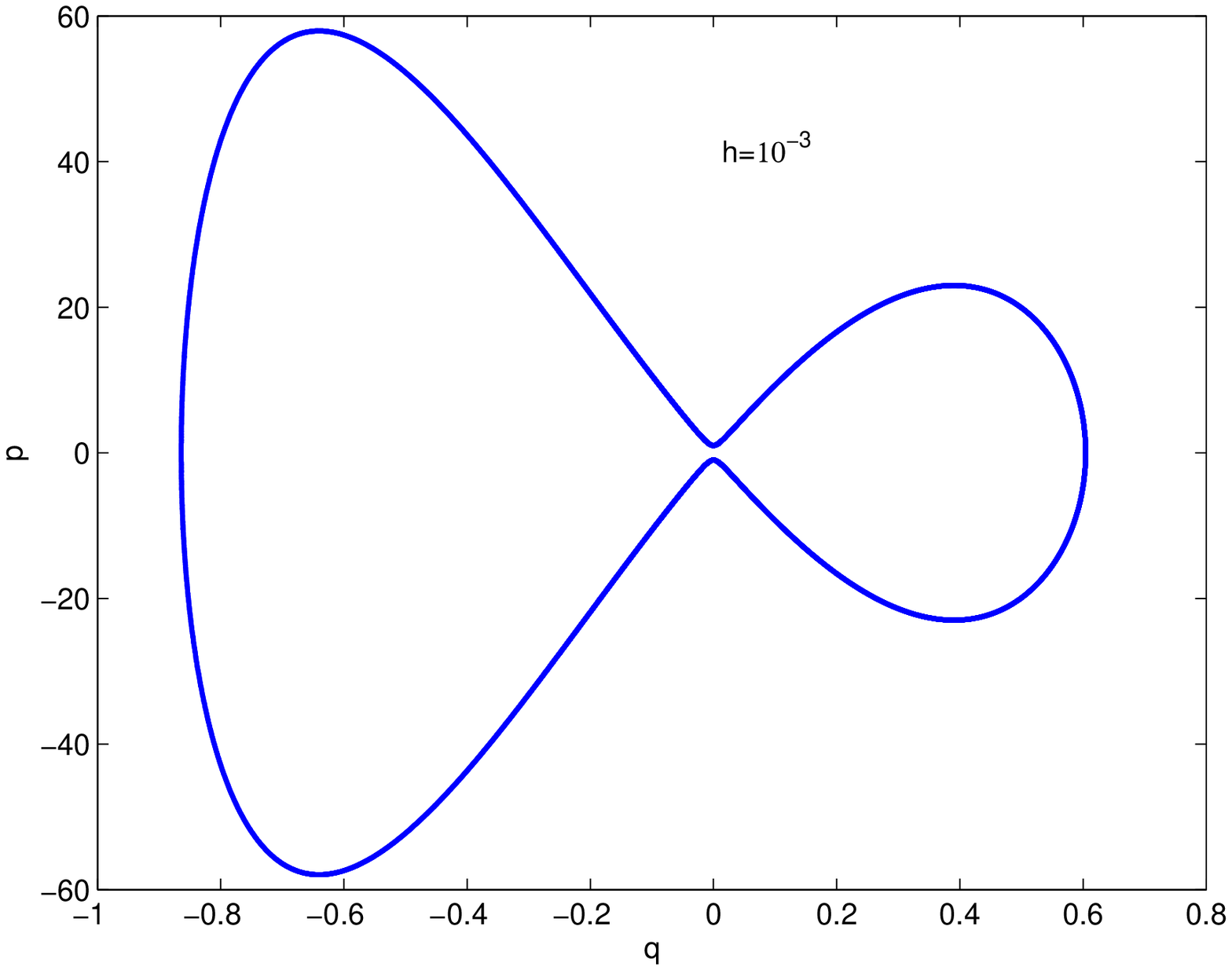}
\qquad
\includegraphics[width=7cm,height=7cm]{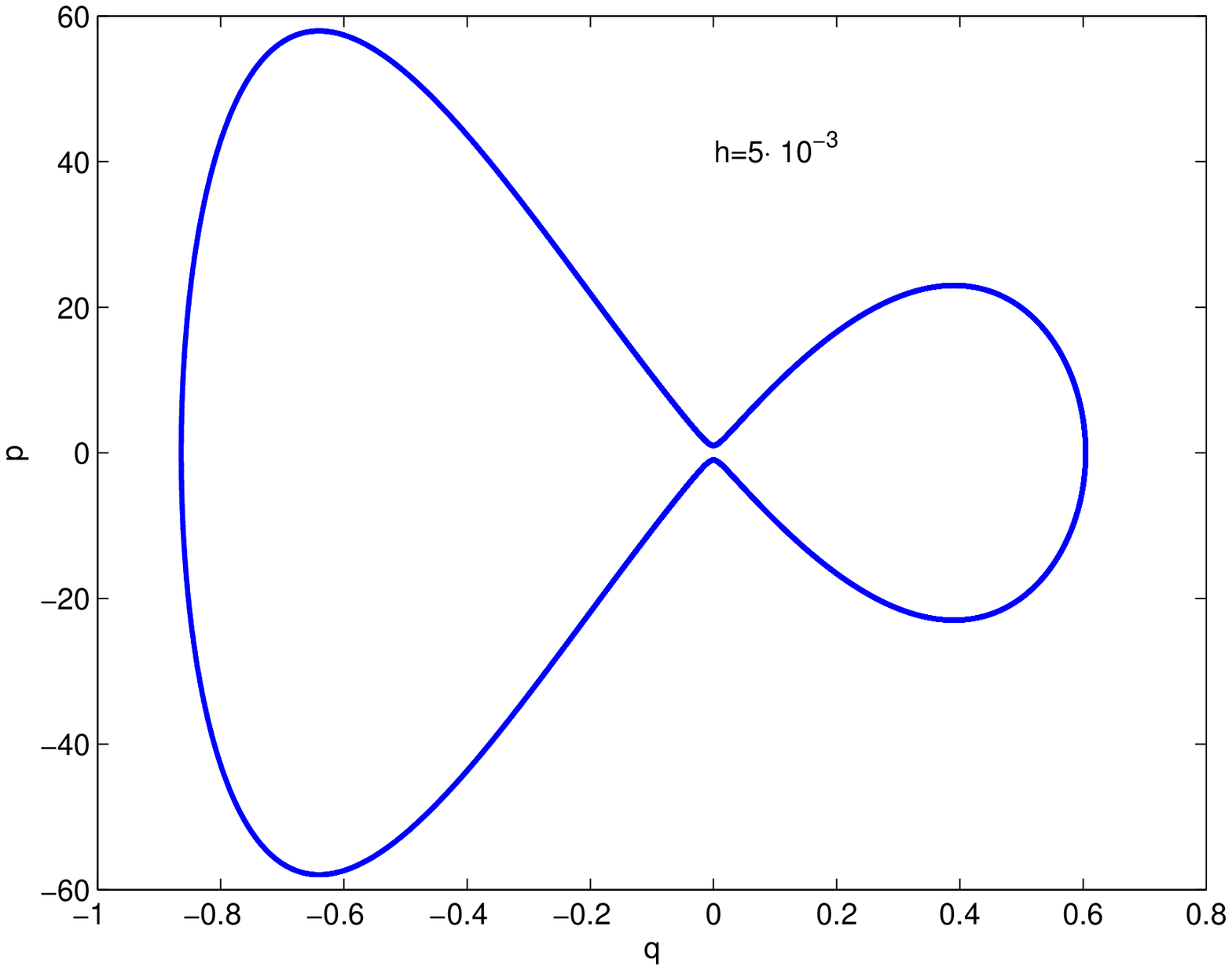}}
\caption{\label{fig1} phase portraits for HBVM(8,2), $h=10^{-3}$
(left) and $h=5\cdot 10^{-3}$ (right).}

\bigskip
\bigskip
\centerline{
\includegraphics[width=7cm,height=7cm]{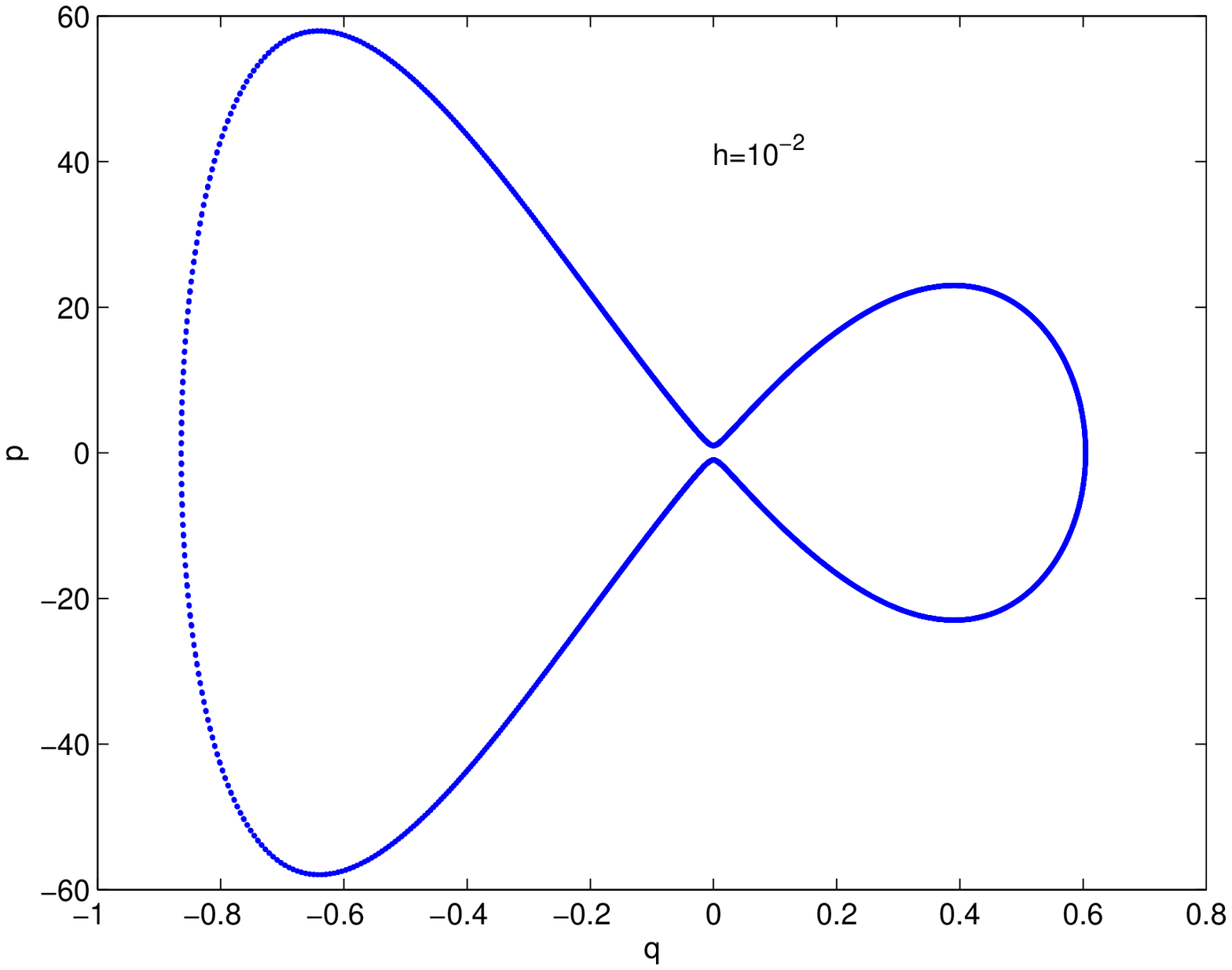}
\qquad
\includegraphics[width=7cm,height=7cm]{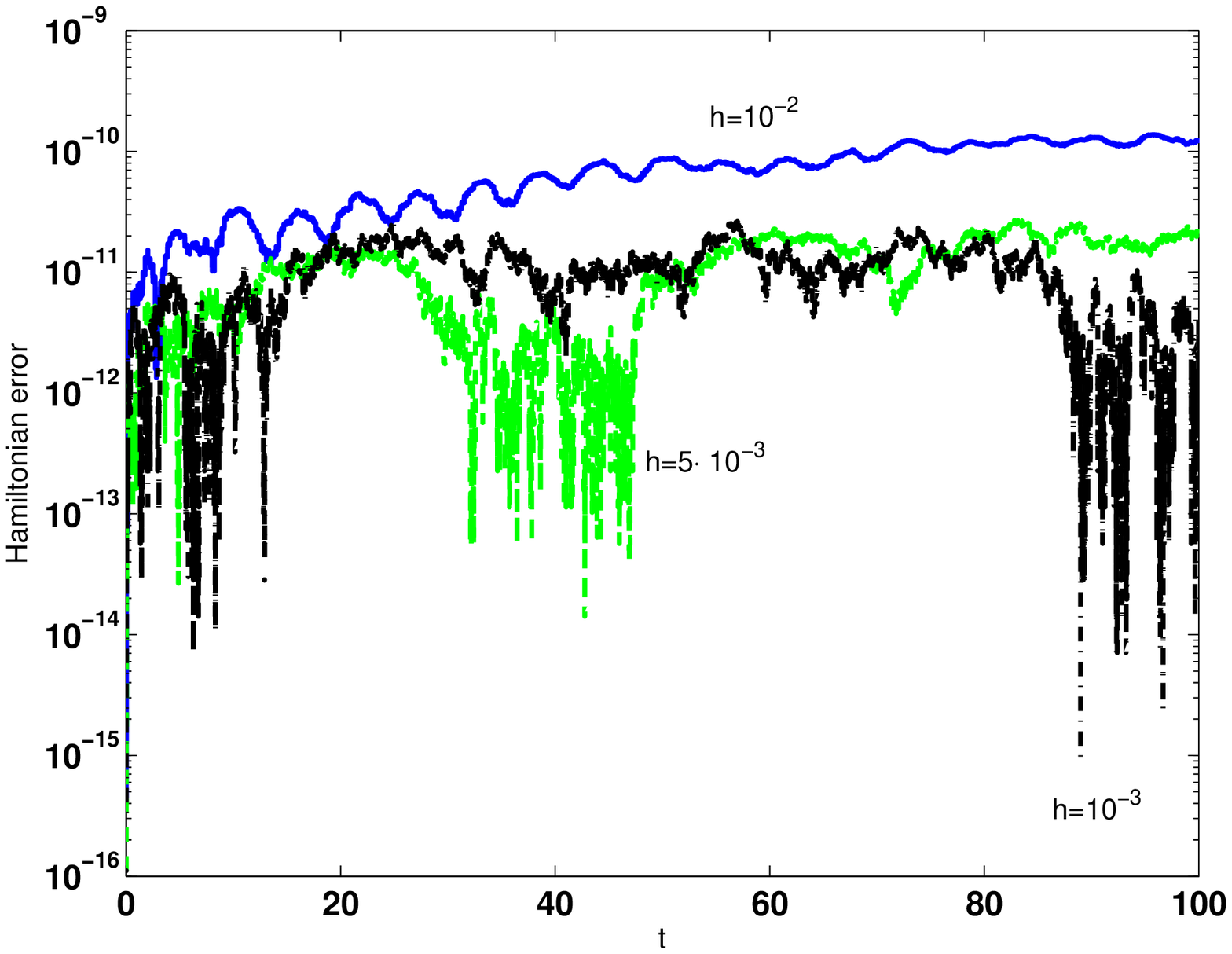}}
\caption{\label{fig2} phase portrait for HBVM(8,2), $h=10^{-2}$
(left) and Hamiltonian error $h=10^{-3}, 5\cdot 10^{-3}, 10^{-2}$
(right).}

\end{figure}

\begin{figure}[ht] \centerline{
\includegraphics[width=7cm,height=7cm]{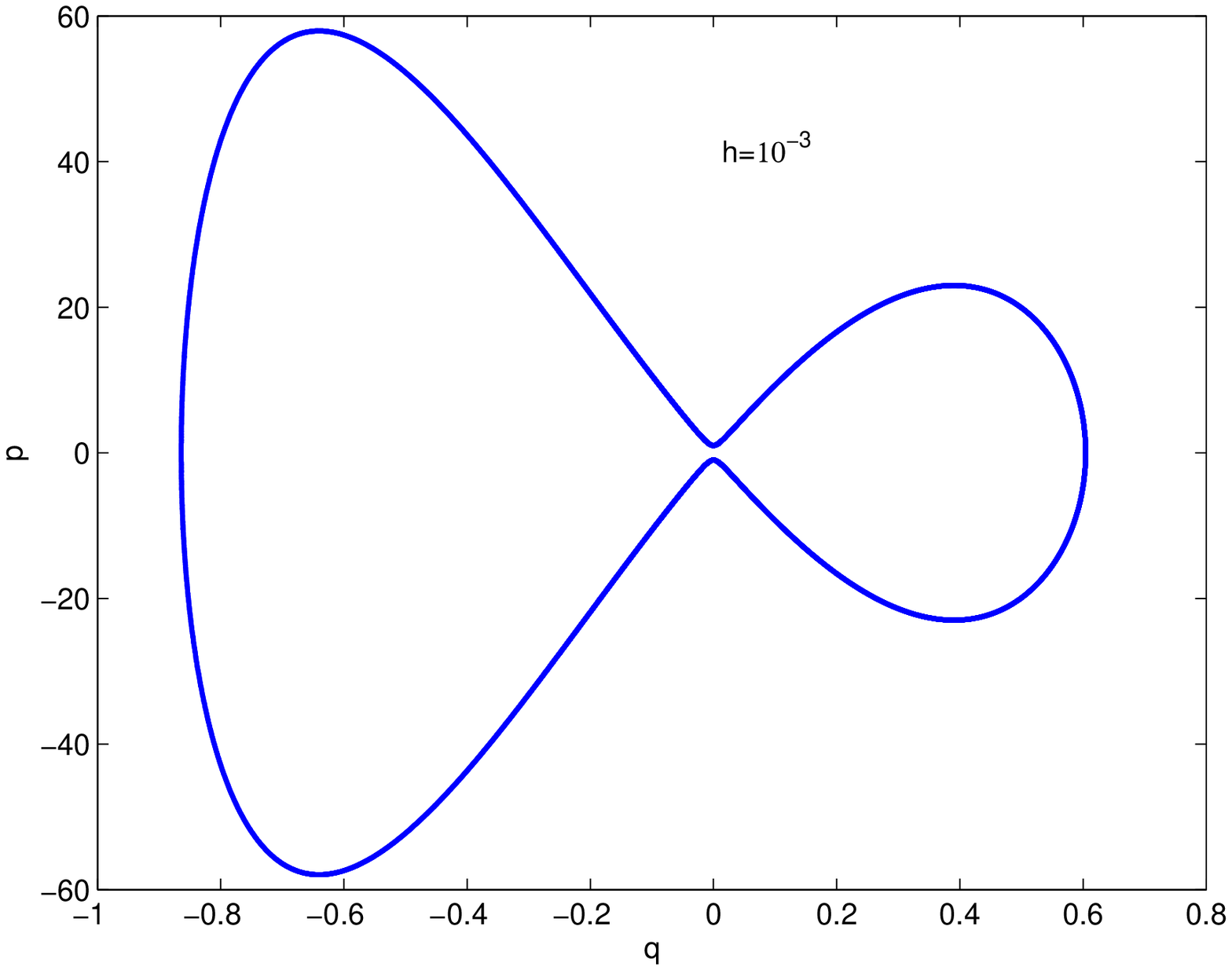}
\qquad
\includegraphics[width=7cm,height=7cm]{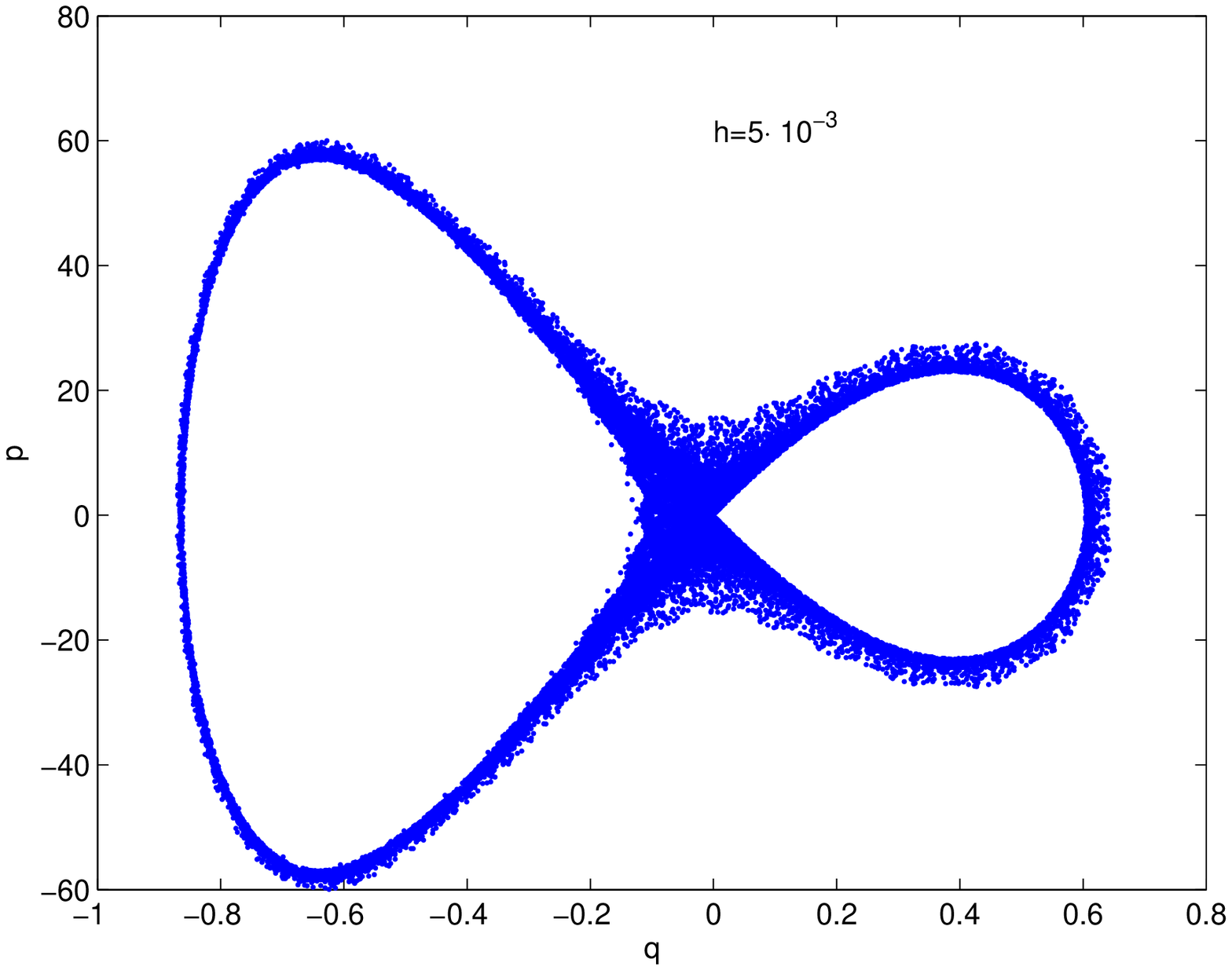}}
\caption{\label{fig3} phase portraits for GAUSS2, $h=10^{-3}$
(left) and $h=5\cdot 10^{-3}$ (right).}

\bigskip
\bigskip
\centerline{
\includegraphics[width=7cm,height=7cm]{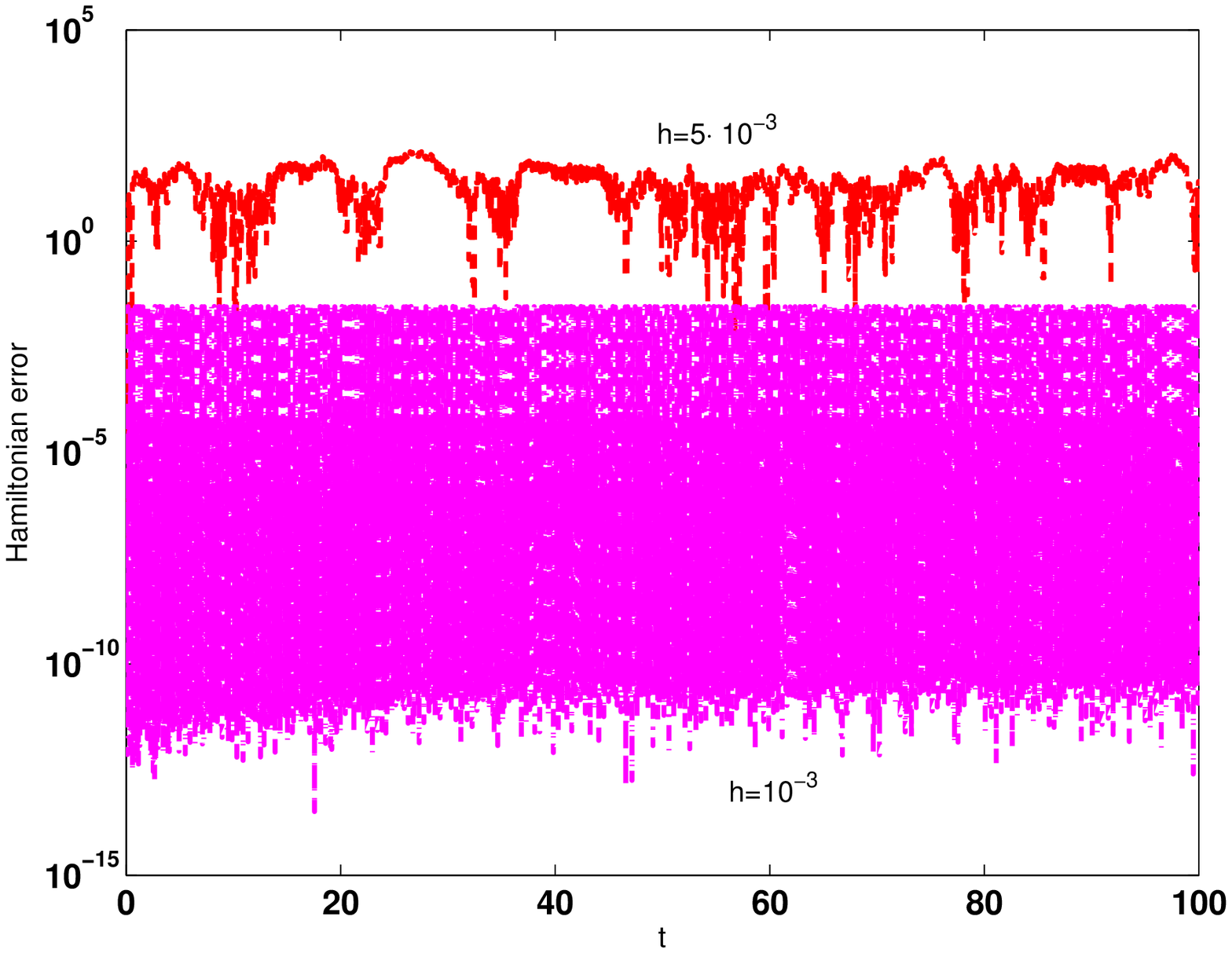}}
\caption{\label{fig4} Hamiltonian error for GAUSS2, $h=10^{-3},
5\cdot 10^{-3}$.}
\end{figure}

Finally, in Figures~\ref{fig1}--\ref{fig4} we plot the phase
portraits of the numerical solutions, for the two methods and the
various stepsizes, along with the corresponding error in the
numerical Hamiltonian. As one can see, the phase portraits of the
HBVM(8,2) method are always correct, whatever the used stepsize
(see Figure~\ref{fig1} and the left plot in Figure~\ref{fig2}), since
the Hamiltonian is conserved (up to round-off errors), as is shown in
the right plot of Figure~\ref{fig2}. This is not true for the
GAUSS2 method, for which the error in the Hamiltonian depends on
the used stepsize, as is shown in Figure~\ref{fig4}, thus causing
drawbacks in the corresponding phase portraits of the numerical
solution, unless the stepsize is very small (see the two plots of
Figure~\ref{fig3}).

From the numerical tests, one can then conclude that the proposed blended
implementation of HBVMs turns out to be robust and efficient. Moreover, the
energy-conserving property of such methods turns out to be very remarkable, with
respect to standard symplectic methods. Finally, the second order formulation of
HBVMs greatly improves their performance.


\end{document}